\documentclass[english]{extarticle}
\usepackage[T1]{fontenc}
\usepackage[latin9]{inputenc}
\usepackage{geometry}
\usepackage{color}
\usepackage{babel}
\usepackage{amssymb}
\usepackage{amsmath}
\usepackage{amsthm}
\usepackage[numbers]{natbib}
\usepackage[unicode=true,pdfusetitle,
 bookmarks=true,bookmarksnumbered=false,bookmarksopen=false,
 breaklinks=false,pdfborder={0 0 1},backref=false,colorlinks=true]
 {hyperref}
\hypersetup{
 linkcolor=blue, urlcolor=marineblue, citecolor=blue, pdfstartview={FitH}, hyperfootnotes=false, unicode=true}

\makeatletter
\numberwithin{equation}{section}
\numberwithin{figure}{section}
\theoremstyle{plain}
\newtheorem{thm}{\protect\theoremname}
\theoremstyle{plain}
\newtheorem{cor}[thm]{\protect\corollaryname}
\newcommand{\DD}{\mathbb{D}}
\usepackage{lettrine} 
\let\myFoot\footnote
\renewcommand{\footnote}[1]{\myFoot{#1\vspace{3mm}}}
\usepackage{lmodern}
\usepackage[T1]{fontenc}
\usepackage{concmath}
\usepackage{ae,aecompl}
\usepackage[T1]{fontenc}
\usepackage{amsmath}
\usepackage{graphicx}

\usepackage{tikz}

\makeatother

\providecommand{\corollaryname}{Corollary}
\providecommand{\theoremname}{Theorem}

\begin{document}
\title{A note on the moments of sequences of complex numbers}

\author{Maher Boudabra, Greg Markowsky\\
Monash University}
\maketitle

\begin{abstract}
We give a short proof that the limsup of the $p$-th root of the modulus of the $p$-th moment of a sequence of complex numbers is equal to the modulus of the maximum of the sequence. This strengthens known results, and provides an analog to a recent result concerning moments of complex polynomials.
\end{abstract}

Given a complex-valued function $f$ defined on the interval $[0,1]$, let the $p$-th moment of $f$ be defined as $M_p = \int_0^1 f(x)^p dx$. In the recent paper \cite{muger2020moments}, the following intriguing result was proved.

\begin{thm} \label{known}
Suppose $f$ is a polynomial with complex coefficients that is not identically 0. Then $\limsup_{p \to \infty}\mid M_{p}\mid^{1/p} > 0$.
\end{thm}

The authors also ventured the natural conjecture that perhaps $\limsup_{p \to \infty}\mid M_{p}\mid^{1/p}$ is equal to $\max_{x \in [0,1]} \mid f(x)\mid$. Note that this theorem and the conjecture may appear simple, given the ease with which they may be proved if $f$ is assumed to be real-valued, but the analytic proof given in \cite{muger2020moments} is far from easy. Upon reflection, this is natural: the moments of $f(x) = e^{2\pi i x}$ are all zero, and so any proof of this theorem must depend heavily upon properties of polynomials. An immediate consequence of this result is that the only polynomial with all moments 0 is the polynomial with all coefficients 0, a result which had earlier been given an algebraic proof (which is also not simple), in \cite{Francoiseetal}.

Motivated by this result, we have considered the analogous result for sequences of complex numbers. Similarly to above, given a sequence of complex numbers $\{z_{n}\}_{n}$ denote by $M_{p}$ for positive integer $p$ the sum $\sum_{n=0}^{+\infty}z_{n}^{p}$. Our result is as follows.

\begin{thm} \label{bigguy}
Suppose a sequence $\{z_{n}\}_{n}$ is in $\ell^q$ for some $q \geq 1$; that is, $\sum_{n=0}^{+\infty}|z_{n}|^{q} < \infty$. Then 

\[
\limsup_{p \to \infty}\mid M_{p}\mid^{\frac{1}{p}}=\max_{n}\mid z_{n}\mid.
\]

\end{thm}

An immediate consequence is the following.

\begin{cor}
Suppose a sequence $\{z_{n}\}_{n}$ is in $\ell^q$ for some $q \geq 1$, and $M_p = 0$ for all $p$. Then $z_n = 0$ for all $n$.
\end{cor}

Before proving our the theorem, several remarks are in order. First, the corollary has already been proven, using a rather technical argument involving C\'es\`aro summation, in \cite{priestley1992complex}. In contrast, our proof of Theorem \ref{bigguy} is very short, and makes use of generating functions and elementary complex analysis. Next we note that the condition that our sequence is in $\ell^q$ is necessary: in \cite{lenard1990nonzero} a remarkable series of nonzero complex numbers is constructed, all of whose moments are 0. Naturally, this sequence is not in $\ell^q$ for any $q$. Finally, we remark that, in our opinion, our theorem lends credence to the aforementioned conjecture of M\"uger and Tuset on moments of polynomials, that $\limsup_{p \to \infty}\mid M_{p}\mid^{1/p} = \max_{x \in [0,1]} \mid f(x)\mid$.

Let us now prove Theorem \ref{bigguy}. Required background for the proof can be found in virtually any complex analysis text (for instance, \cite{Rudin2001} or \cite{spiegel2015complex}). Let us begin by assuming that $\{z_{n}\}_{n}$ is in $\ell^1$, and by multiplying by a constant we may also assume that $\max_{n}\mid z_{n}\mid = 1$ (we are assuming that $z_n \neq 0$ for at least one $n$, otherwise the result is trivial). For $w\in\mathbb{D}$, define

\begin{equation} \label{fdefn}
    f(w) = \sum_{n=1}^\infty \frac{z_n w}{1-z_n w}.
\end{equation}

Since $|1-z_n w|$ is bounded uniformly away from $0$ on compact subsets of $\DD$, the fact that $\{z_{n}\}_{n} \in \ell^1$ implies that the series in (\ref{fdefn}) converges locally uniformly, and $f$ is therefore analytic on $\DD$. By expanding each term in the series into a geometric series and swapping the order of summation (which is straightforward to justify using the locally uniform convergence), we find that the power series expansion for $f$ is simply the generating function for the sequence $M_p$, namely

\begin{equation} \label{fpower}
    f(w) = \sum_{p=1}^\infty M_p w^p.
\end{equation}

However, it is evident from (\ref{fdefn}) that $\lim_{w \to 1/z_n} |f(w)| = \infty$ if $|z_n|=1$, and we have assumed that there is at least one such $z_n$. This implies that $f$ can not be extended to an analytic function on a disk centered at the origin of radius greater than one, and therefore that the radius of convergence of the power series in (\ref{fpower}) is 1. In other words,

\[
\limsup_{p \to \infty}\mid M_{p}\mid^{\frac{1}{p}}=1=\max_{n}\mid z_{n}\mid.
\]

Now we suppose that $\{z_{n}\}_{n}$ is not in $\ell^1$ but is in $\ell^q$ for some positive integer $q>1$, and we still assume $\max_{n}\mid z_{n}\mid = 1$. Let $v_n = z_n^q$, and note that $M_p(v) = M_{pq}(z)$. The prior argument applies to $\{v_{n}\}_{n}$, since this sequence is in $\ell^1$, and so

\begin{equation}
    \limsup_{p\to \infty} \mid M_{p}(z)\mid^{\frac{1}{p}} \geq \limsup_{p\to \infty}\mid M_{pq}(z)\mid^{\frac{1}{pq}} = \Big(\limsup_{p\to \infty} \mid M_{p}(v)\mid^{\frac{1}{p}}\Big)^{\frac{1}{q}} = 1.
\end{equation}

On the other hand, the quantities $\mid M_{p}(z)\mid$ are uniformly bounded for large enough $p$ (since $\max_{n}\mid z_{n}\mid = 1$ and $\{z_{n}\}_{n} \in \ell^q$), so that $\limsup_{p\to \infty} \mid M_{p}(z)\mid^{\frac{1}{p}} \leq 1$. The result follows.

\bibliographystyle{plain}
\bibliography{Maheref}

\begin{thebibliography}{1}

\bibitem{Francoiseetal}
J.P. Francoise, F.~Pakovich, Y.~Yomdin, and W.~Zhao.
\newblock Moment vanishing problem and positivity: Some examples.
\newblock {\em Bulletin des Sciences Math\'ematiques}, 135(1):10--32, 2011.

\bibitem{lenard1990nonzero}
A.~Lenard.
\newblock A nonzero complex sequence with vanishing power-sums.
\newblock {\em Proceedings of the American Mathematical Society},
  108(4):951--953, 1990.

\bibitem{muger2020moments}
M.~M{\"u}ger and L.~Tuset.
\newblock On the moments of a polynomial in one variable.
\newblock {\em Indagationes Mathematicae}, 31(1):147--151, 2020.

\bibitem{priestley1992complex}
W.~Priestley.
\newblock Complex sequences whose "moments" all vanish.
\newblock {\em Proceedings of the American Mathematical Society},
  116(2):437--444, 1992.

\bibitem{Rudin2001}
W.~Rudin.
\newblock {\em Real and complex analysis (3rd ed.)}.
\newblock McGraw-Hill Education, 2001.

\bibitem{spiegel2015complex}
M.~Spiegel, S.~Lipschutz, J.~Schiller, and D.~Spellman.
\newblock {\em Complex variables: With an introduction to conformal mapping and
  its applications}.
\newblock McGraw-Hill, 2015.

\end{thebibliography}

\end{document}